# Oriented Convex Containers of Polygons

R. Nandakumar,
Amrita School of Arts and Sciences,
Idapalli North, Kochi 682024, India
nandacumar@gmail.com

**Abstract:** We consider the optimal containment of polygonal regions within convex containers of with the special property of 'oriented-ness' and derive preliminary results.

**Introduction:** We define an oriented planar region as a convex region that enables us to consistently specify a preferred direction on the plane. Examples are isosceles triangles (the natural choice of orientation is the bisector of its apex), the rectangle (the orientation can be chosen to be given by the longer side) and the ellipse (the major axis as orientation). There could also be oriented pentagons, hexagons and so forth.

Here, we study containment of any given convex region *optimally* with various types of oriented regions. For each type of container, we try to optimize the container in two different ways – to minimize the area of the container and to minimize the perimeter. We consider rectangles, isosceles triangles and ellipses as containers.

This class of problems illustrate that the area and perimeter of convex regions are not totally unrelated quantities. Indeed, for a general input region and any specified type of oriented container for this region, we find that the container with least area (perimeter) is often the container with least perimeter (area) as well. Even in cases where an oriented container achieves a minimum for only one of the two quantities, we would often (although not always) be quite close to the minimum of the other quantity as well. In what follows, we make these statements more precise by calculating, given any convex region to be contained, how different the orientations of the containers minimizing the two quantities are.

1. *Rectangle containers of least area and perimeter for a given convex region*

For a convex polygon, let us call its rectangular container of minimum area $R_A$ and the rectangular container of minimum perimeter $R_P$. We note an important property:

**Lemma 1:** Given any convex region C, its $R_A$ necessarily has an edge that overlaps at least one edge of C. This property holds for its $R_P$ as well. The edges of C that overlaps with its two optimal rectangle containers need not be the same.

Note: Lemma 1 can be proved easily with basic trigonometry and is the basis of a rotating calipers approach to find the least area and least perimeter bounding rectangles of any set of points([1]). Here we are not interested in finding these rectangles efficiently but in finding how they could vary in orientation for the same input convex region.

We first describe a triangle for which $R_P$ and $R_A$ differ in orientation by almost 45 degrees.

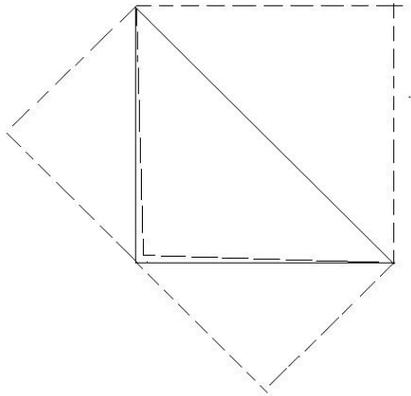

Figure 1

For an isosceles right triangle – call it T - with unit side (drawn with solid lines above), $R_P$ is a square of side 1 shown above. There are two differently oriented rectangles, as shown above, both of area 1 unit that could be the $R_A$ for the isosceles right triangle. Let us call the square R1 and the rectangle tilted at 45 degrees, R2. R1 has area 1 and perimeter 4. R2 has area 1 and perimeter nearly 4.242.

Now, deform T slightly as shown with dashed lines above, keeping it isosceles but with the apex angle slightly more than 90 degrees. Call this new triangle T'. It is now easy to see that by parallel displacing a long side of R2 slightly, we have a rectangle, say R2', that just contains T'. A rectangle R1' in the neighborhood of square R1 (very slightly tilted with respect to R1) also contains T'. But we can easily see that R2' has slightly less area than R1' and so it is the sole $R_A$ of T'. But R1', which is only a slight deformation of the unit square, has perimeter clearly less than R2' and is the unique $R_P$ of T'. So, T' is a triangle that has $R_A$ and $R_P$ differing in orientation by almost 45 degrees.

Moving to more general input regions, we show another region for which $R_A$ and $R_P$ have very different orientations. We follow a construction suggested by Friedman [4].

Consider a square. From a pair of opposite corners we keep 'shaving' off isosceles triangles to leave a progressively thinner hexagon aligned along one of the diagonals of the square (this diagonal is the diameter of the hexagon). See the following figure that shows an intermediate stage in this thinning process.

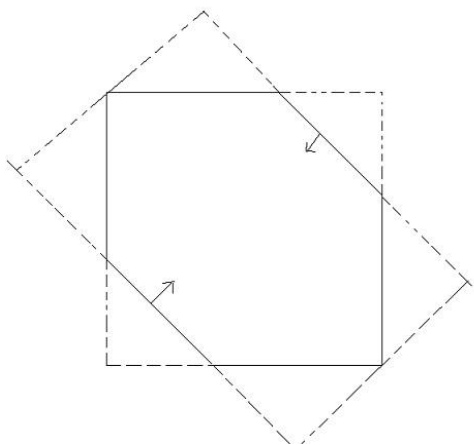

Figure 2

In the initial stages of this continuous transformation, both the least area and least perimeter rectangles that contain the hexagon are the original square itself. It is also evident that when the hexagon gets very thin, both the least area and least perimeter rectangular containers are aligned along the diameter of the hexagon (both these rectangular containers are shown above). It can be checked numerically that for $R_A$ and $R_P$, the switch from the full square to a rectangle aligned along the hexagon diameter happens at different stages of the thinning. So between these two switch points, the two rectangular containers differ in orientation by 45 degrees.

**Question:** Are there convex regions for which $R_A$ and $R_P$ vary in orientation by more than 45 degrees – say even by an angle arbitrarily close to 90 degrees?

Answer: The answer is yes.

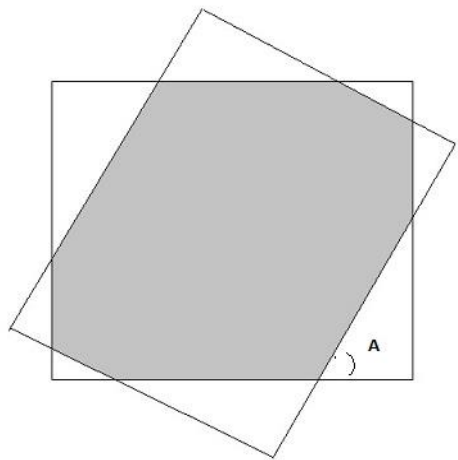

Figure 3

In the figure above are shown two rectangles R1 (axis parallel) and R2 (tilted) with the following properties:

 - The perimeter of R1 < perimeter of R2.

- The area of R1 > area of R2

- The diagonal of R1 > length of R2

- R2 is oriented at an angle A to T1 and A is more than 45 degrees and the centers of both are coincident.

Consider the shaded octagonal region which is the intersection of R1 and R2 – call this O. Since its $R_P$ and $R_A$ have to share at least an edge with it we see readily that R1 is the $R_P$ and R2 is the $R_A$ of O. And the difference in orientation between them is more than 45 degrees. Clearly, A can have a range of values greater than 45 degrees for which such an octagon exists. Indeed, we find that one can choose R1 and R2 in such a way that for A arbitrarily close to 90 degrees, we can form such an octagonal region.

A specific example: R1 has dimensions 10 X 9.9; R2 has dimensions 11 X 8.99. R1 has greater area and R3 has greater perimeter. For these rectangles, for values of angle A up to almost 83 degrees, we can form such an octagon.

2. *Isosceles triangle containers of least area and least perimeter*

The most basic question here: Given a general triangle T, to find isosceles triangles of least area that contains T (we call this isosceles triangle $T_A$) and the least perimeter isosceles triangle that contains T (call this isosceles triangle $T_P$). We had used the following claim in the earlier version of this paper.

**Claim:** Given any triangle T, both its optimal isosceles triangle containers $T_A$ and $T_P$ necessarily share an angle with T.

It has since been clarified that the claim does not hold in general if T is obtuse. So, as of now, we do not have an answer to how much the orientations of the two optimal isosceles triangles of a given T can vary.

**Further Question:** It could also be of interest, given some specified angle α, to construct general convex polygons for which the difference in orientation of $T_A$ and $T_P$ is provably α.

**Remark:** While for any convex polygonal region C, at least one edge of $R_A$ (and $R_P$) lies flush over an edge of C (Lemma 1), we cannot make a similar assertion for the sides of the isosceles triangles $T_A$ and $T_P$. our numerical experiments show that neither $T_A$ nor $T_P$ necessarily has a side flush with C. Due to this, it appears that we cannot readily proceed as in figure 3 and take the intersection of two isosceles triangles with suitable areas and perimeters and with widely varying orientations to form a convex region for which the $T_A$ and $T_P$ provably have widely varying orientations.

**Right Triangle Containers of least area and perimeter:** Right triangles, though in general lacking a direction of symmetry, can be considered oriented – with any right triangle a rectangle is naturally associated and we could take the orientation of this rectangle. We now consider the least area and least perimeter right triangles that contain a given convex region (call these triangles $RT_A$ and $RT_P$ respectively) and look for convex regions for which these optimal containers differ most in orientation.

We revisit the construction of figure 3 for finding convex regions for which the two rectangular containers have widely varying orientations. See figure below.

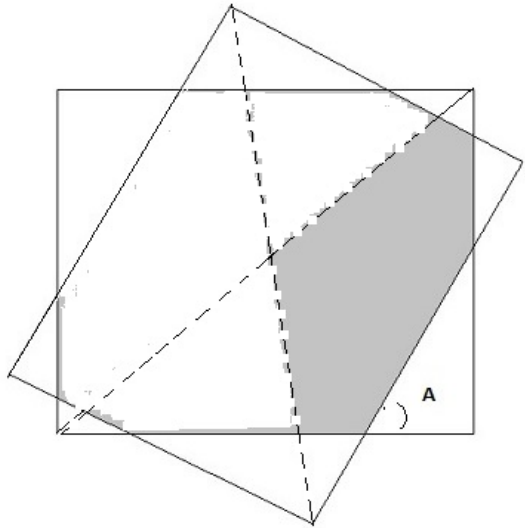

Figure 7

We recall that in figure 3, the horizontal rectangle has less perimeter and greater area than the tilted rectangle. So in figure 5, the right triangle that is half of the horizontal rectangle has less perimeter and more area than the right triangle that is half of the tilted rectangle (the hypotenuses of both right triangles are shown as dashed lines). We guess that for the shaded hexagonal region above, the one from two right triangles with less perimeter is the $RT_P$ and the one with less area is the $RT_A$ and hence further guess that we can thus construct convex regions for which the two optimal right triangle containers can have widely varying orientations.

### 3. *Elliptical Containers of least area and perimeter*

Here we look for ellipses of least area and least perimeter ($E_A$ and $E_P$) that contain an input convex region. We numerically calculated $E_A$ and $E_P$ for triangles and quadrilaterals as input. As mentioned above, we consider the direction of the major axis of the ellipse to be its orientation.

Note: There is no known closed form for the perimeter of an ellipse so we use a close numerical approximation. If a and b are the major and minor axes of an ellipse, its perimeter is given approximately by the Gauss Kummer series.

$$p = \pi (a+b) ( 1 + h/4 + h^2/64 + h^3/256 + 25\, h^4/16384 + ….) \qquad [3]$$

where $h = (a-b)^2 / (a+b)^2$. In our calculations, we stop with the quartic term. The lack of a closed form expression for perimeter of ellipse does not necessarily mean that the perimeters of two ellipses cannot be compared.

**Summary of Findings:**

- For most shapes of triangles as input, $E_A$ and $E_P$ are NOT identical (Let us recall: $T_A$ and $T_P$ are identical for most triangles; $R_A$ and $R_P$ too are identical for most triangles). The centers of $E_A$ and $E_P$ lie considerably apart for many triangles. But the orientations of $E_A$ and $E_P$ are found to be very close for any general triangle – indeed, we find that they are never apart by more than 2 degrees. Further, we note that the orientations of these elliptical containers differ from that of the longest side of T by only a few degrees.

- For any convex quadrilateral Q, the orientations of $E_A$ and $E_P$ are again close – the maximum difference found experimentally is less 12 degrees (for most quadrilaterals, this difference is much less although finite). The orientations of $E_A$ and $E_P$ of any quadrilateral differ from that of the diameter of the quadrilateral by an angle within 30 degrees (for most quadrilaterals, this difference is within 15 degrees).

- For parallelograms as input, $E_A$ and $E_P$ usually have different orientations but as we checked numerically, parallelograms are not the quadrilaterals for which the orientations of $E_A$ and $E_P$ are maximally different.

**Further Questions:** Are there convex regions for which the difference in orientation between the two containers $E_A$ and $E_P$ is provably large, say 45 degrees or higher? We don't have an answer as of now; we don't know if the construction shown in figure 3 can be applied with ellipses replacing rectangles. Indeed we suspect that unlike in the case of rectangular containers, there may be a finite upper bound on the difference in orientations of the two optimal elliptical containers. If the suspicion is correct, one can also ask if the ellipse is the oriented container for which the upper bound of the difference in orientation of the two types of container is the least.

**Acknowledgements:** Thanks to K Sheshadri for his very helpful comments. Thanks to Pinaki Majumdar and HRI, Allahabad for their generous support.